\documentclass[11pt,a4paper,fleqn]{article}

\usepackage[fleqn]{amsmath}
\usepackage{amsthm}
\usepackage{graphicx}
\usepackage{url}
\usepackage{enge}
\usepackage{fancyvrb}

\let\epsilon=\varepsilon
\DeclareMathOperator{\poly}{poly}
\providecommand{\binom}[2]{{#1 \choose #2}}

\providecommand{\from}{\colon}
\providecommand{\suchthat}{:}

\providecommand{\dunion}{%
  \mathop{\ooalign{$\cup$\crcr\hfil\raise0.7ex\hbox{.}\hfil\crcr}}}

\newtheorem{theorem}{Theorem}
\newtheorem{lemma}{Lemma}

\theoremstyle{definition}
\newtheorem{claim}{Claim}

\title{Bipartite Multigraphs with Expander-Like Properties}

\author{Lars Engebretsen}

\institute{%
  Department of Numerical Analysis and Computer Science\\
  Royal Institute of Technology\\
  SE-100 44 Stockholm\\
  SWEDEN\\
  E-mail: enge@kth.se}

\date{December 2004}

\begin{document}

\maketitle

\section{Introduction}

A graph with vertex set~$V$ and edge set~$E$ is called a
\emph{$(d,c)$-expander} if the maximum degree of a vertex is~$d$ and,
for every set $W \subset V$ of cardinality at most $|V|/2$, the
inequality $|\{w \in V \setminus W \suchthat \{v,w\} \in E\}| \geq
c|V|$ holds. This note considers a related combinatorial question:
\begin{quote}
  For which integers~$d$ and functions~$f_d$ does there exist, for
  every large enough~$v$, a bipartite $d$-regular multigraph on $2v$
  nodes with node sets $V$ and~$W$ having the following property: For
  every $U \subseteq V$ and every $U \subseteq W$, the cardinality of
  the set of neighbours of~$U$ is at least $f_d(|U|)$?
\end{quote}
Graphs with the above property seem to behave well also with respect
to other, more complicated, expansion-type properties. Indeed, the
author was motivated to study this question by a paper communicated to
him in May 2002 (the latest version of the paper is available from URL
\url{http://www-math.mit.edu/~vempala/papers/tspinapprox.ps}). In this
paper, Papadimitriou and Vempala established approximation hardness of
TSP with triangle inequality using as a tool in their construction the
fact that for $d = 6$ and
\begin{displaymath}
  f_6(u) = \begin{cases}
    2u        & \mbox{$0   \leq u \leq v/4$,} \\
    u + v/4   & \mbox{$v/4 \leq u \leq v/2$,} \\
    u/2 + v/2 & \mbox{$v/2 \leq u \leq v$,}
  \end{cases}
\end{displaymath}
there exist bipartite multigraphs with the properties described in
the above question. In this paper, we prove the following theorem:
\begin{theorem}
  \label{thm:main}
  For $d \in \{5,6,7,8\}$ and functions $f_d$ as described below,
  there exists, for every large enough~$v$, a bipartite $d$-regular
  multigraph on $2v$ nodes with node sets $V$ and~$W$ having the
  property that for every $U \subseteq V$ and every $U \subseteq W$
  the cardinality of the set of neighbours of~$U$ is at least
  $f_d(|U|)$.
  \begin{gather}
    \label{eq:f5}
    f_5(u) = \begin{cases}
      2u            & \mbox{$0 \leq u \leq 3v/20$,} \\
      4u/3 + v/10   & \mbox{$3v/20 \leq u \leq 3v/10$,} \\
      u + v/5       & \mbox{$3v/10 \leq u \leq v/2$,} \\
      3u/4 + 13v/40 & \mbox{$v/2 \leq u \leq 7v/10$,} \\
      u/2 + v/2     & \mbox{$7v/10 \leq u \leq v$,}
    \end{cases}
    \\
    \label{eq:f6}
    f_6(u) = \begin{cases}
      5u/2         & \mbox{$0 \leq u \leq v/10$,} \\
      5u/3 + v/12  & \mbox{$v/10 \leq u \leq v/4$,} \\
      u + v/4      & \mbox{$v/4 \leq u \leq v/2$,} \\
      3u/5 + 9v/20 & \mbox{$v/2 \leq u \leq 3v/4$,} \\
      2u/5 + 3v/5  & \mbox{$3v/4 \leq u \leq v$,}
    \end{cases}
    \\
    \label{eq:f7}
    f_7(u) = \begin{cases}
      3u               & \mbox{$0 \leq u \leq v/10$,} \\
      2u + v/10        & \mbox{$v/10 \leq u \leq 3v/20$,} \\
      21u/15 + 19v/100 & \mbox{$3v/20 \leq u \leq 3v/10$,} \\
      u + 31v/100      & \mbox{$30v/100 \leq u \leq 39v/100$,} \\
      15u/21 + 59v/140 & \mbox{$39v/100 \leq u \leq 3v/5$,} \\
      u/2 + 11v/20     & \mbox{$3v/5 \leq u \leq 7v/10$,} \\
      u/3 + 2v/3       & \mbox{$7v/10 \leq u \leq v$.}
    \end{cases}
    \\
    \label{eq:f8}
    f_8(u) = \begin{cases}
      3u           & \mbox{$0 \leq u \leq v/10$,} \\
      2u+v/10      & \mbox{$v/10 \leq u \leq v/5$,} \\
      5u/4 + v/4   & \mbox{$v/5 \leq u \leq v/3$,} \\
      4u/5 + 2v/5  & \mbox{$v/3 \leq u \leq v/2$,} \\
      u/2 + 11v/20 & \mbox{$v/2 \leq u \leq 7v/10$.} \\
      u/3 + 2v/3   & \mbox{$7v/10 \leq u \leq v$.}
    \end{cases}
  \end{gather}
\end{theorem}
The proof of the theorem relies partly on claims that have been
verified by a computer. Hence, a second contribution in this paper is
the description of a fairly general methodology for devising
computer-assisted proofs for a wide class of mathematical claims.

\subsection{The probabilistic method}

The probabilistic method~\cite{AloSpe91} is particularly well-suited
for proving that there exists, in some large class of objects, an
object with certain specified properties. Typically, the argument
proceeds by first selecting a random object from the class and then
estimating the probability that the selected object \emph{does not
  have} the sought properties. If this probability can be shown to be
strictly less than one, the probability that the selected object
\emph{has} the property is strictly positive and hence \emph{there
  exists} an object with the sought properties.

The method is highly non-constructive. For our case, we in fact have
no idea whatsoever how to actually construct a graph with the sought
properties in time polynomial in the number of vertices in the graph.

\subsection{Computer-assisted proofs}

As mentioned in the previous section, a critical component in the
probabilistic method is to show that the probability of some event is
strictly less than one. In our case, the function expressing this
probability is fairly complicated, although continuously differentiable
almost everywhere. In principle, straightforward but tedious analysis
of first and second order derivatives of the function could be used to
prove that it is strictly less than one. We feel, however, that the
contribution to the community from such a proof is very minor. Instead
we resort to a computer-assisted proof and argue that this method of
proof should be accepted in cases like ours.

First, what is a computer assisted proof of some statement ``$f(x) <
c$ for all $x \in [a,b]$''? Simply evaluating the function at some
points is clearly not enough---the function may assume other,
dangerous, values at the points where it was not evaluated. Given some
bound, proven by a conventional mathematical proof, of the form
``$|f'(x)| < C$ for all $x \in [a,b]$'', we could argue that it is
enough to evaluate the function at points that are sufficiently close
since a Taylor expansion then bounds the value of the function at all
points. However, floating point computations done by computers are not
accurate, and the latter argument above fails to take into account
possible influences of round-off errors.

A solution to the problem of round-off errors is to use interval
arithmetic~\cite{Moo66}. The main idea behind computations with
interval arithmetic is to compute not with single numbers but rather
with intervals. When some function~$f$ is applied to some
interval~$I$, the result is an interval that contains $f(x)$ for every
$x \in I$. Hence, interval arithmetic is particularly well suited for
verifying claims that are of the form ``$f(x) < c$ for $x \in I$''. If
the interval~$I$ is large, then the result of computing $f(I)$ is
usually also a large interval; in particular, $I$ could also
contain~$c$. The solution to this problem is to split $I$ into
sub-intervals that are sufficiently small and then compute $f(I_k)$
for each sub-interval~$I_k$. The computer-assisted part of our proofs
works precisely in this way and was inspired by Uri Zwick's work on
optimal approximation algorithms for certain constraint satisfaction
problems~\cite{Zwi02}.

We argue that computer-assisted proofs should be accepted in cases
like ours. To examine and ascertain the correctness of the program
presented in this paper is the same thing as examining a conventional
mathematical proof. It is, of course, true that the correctness of the
computer-assisted \emph{verification} of the claims that the program
verifies, requires the assumption that the program is correctly
compiled and that all used library routines are correctly written. It
is in principle possible to construct a compiler that automatically
includes logging facilities in the executable program; the output from
such a system could then in principle be verified step by step to
ascertain the correctness of the computation.

\section{Proof of Theorem~\ref{thm:main}}

We select a $d$-regular bipartite multigraph on $2v$~vertices by
selecting one perfect matching in a bipartite graph on $dv + dv$
vertices uniformly at random. From this perfect matching, the
$d$-regular bipartite graph is constructed by identifying groups of
$d$~vertices in the ``big'' bipartite graph with single vertices in
the sought $d$-regular bipartite graph.

The analysis of the construction proceeds by estimating the
probability that such a randomly chosen graph does not have the
desired properties, i.e., the probability that there is some $U
\subseteq V$ or some $U \subseteq W$ such that the cardinality of the
set of neighbours of~$U$ is no more than $f_d(|U|)$.

Fix a set $U \subseteq V$ of size~$u$ and a set $N \subseteq W$ of
size~$n$. When $u > n$, the pigeon-hole principle implies that
$U$~cannot have neighbours only in~$N$. When $u \leq n$, the
probability that $U$ has neighbours only in~$N$ is
\begin{displaymath}
  \binom{dn}{du} \bigg/ \binom{dv}{du}
\end{displaymath}
since there are in total $\binom{dv}{du}$ ways to choose the
neighbours of~$U$ and $\binom{dn}{du}$ of those choices result in
neighbours only in~$N$. It turns out that we need different methods to
bound this probability depending on how close $(u,n)$ is to the
``extreme points'' of a certain region. Therefore, we define
\begin{gather*}
  \Omega = \{(u,n) \suchthat u \leq n \leq f_d(u)\}
  \mbox{,}\\
  \Omega_\delta = \{(u,n) \in \Omega \suchthat \min\{u,v-n\} \leq \delta v\}
  \mbox{,}\\
  \Omega' = \Omega \setminus \Omega_{\delta}
  \mbox{.}
\end{gather*}
With this notation, the probability that there is some $U \subseteq V$
or some $U \subseteq W$ such that the cardinality of the set of
neighbours of~$U$ is no more than $f_d(|U|)$ can be written
\begin{displaymath}
  2
  \sum_{(u,n)\in\Omega}
  \sum_{\substack{U \in V\\|U|=u}}
  \sum_{\substack{N \in W\\|N|=n}}
  \binom{dn}{du} \bigg/ \binom{dv}{du}
\end{displaymath}
where the factor~$2$ above comes from the fact the we consider not
only neighbour sets of $U \subseteq V$ but also neighbour sets of $U
\subseteq W$. Let
\begin{equation}
  \label{eq:P}
  P(u,n)
  = 
  \binom{v}{u} \binom{v}{n} \binom{dn}{du} \bigg/ \binom{dv}{du}
  = 
  \binom{v}{u} \binom{v}{n} \frac{(dn)!(dv-du)!}{(dn-du)!(dv)!}
  \,\mbox{.}
\end{equation}
Using $P$ and the $\Omega$'s, the probability that there is some $U
\subseteq V$ or some $U \subseteq W$ such that the cardinality of the
set of neighbours of~$U$ is no more than $f_d(|U|)$ can be upper
bounded by $2|\Omega_{\delta}| \max_{(u,n)\in\Omega_{\delta}} P(u,n) +
2|\Omega'| \max_{(u,n)\in\Omega'} P(u,n)$.

\begin{figure}[p]
  \centerline{\includegraphics{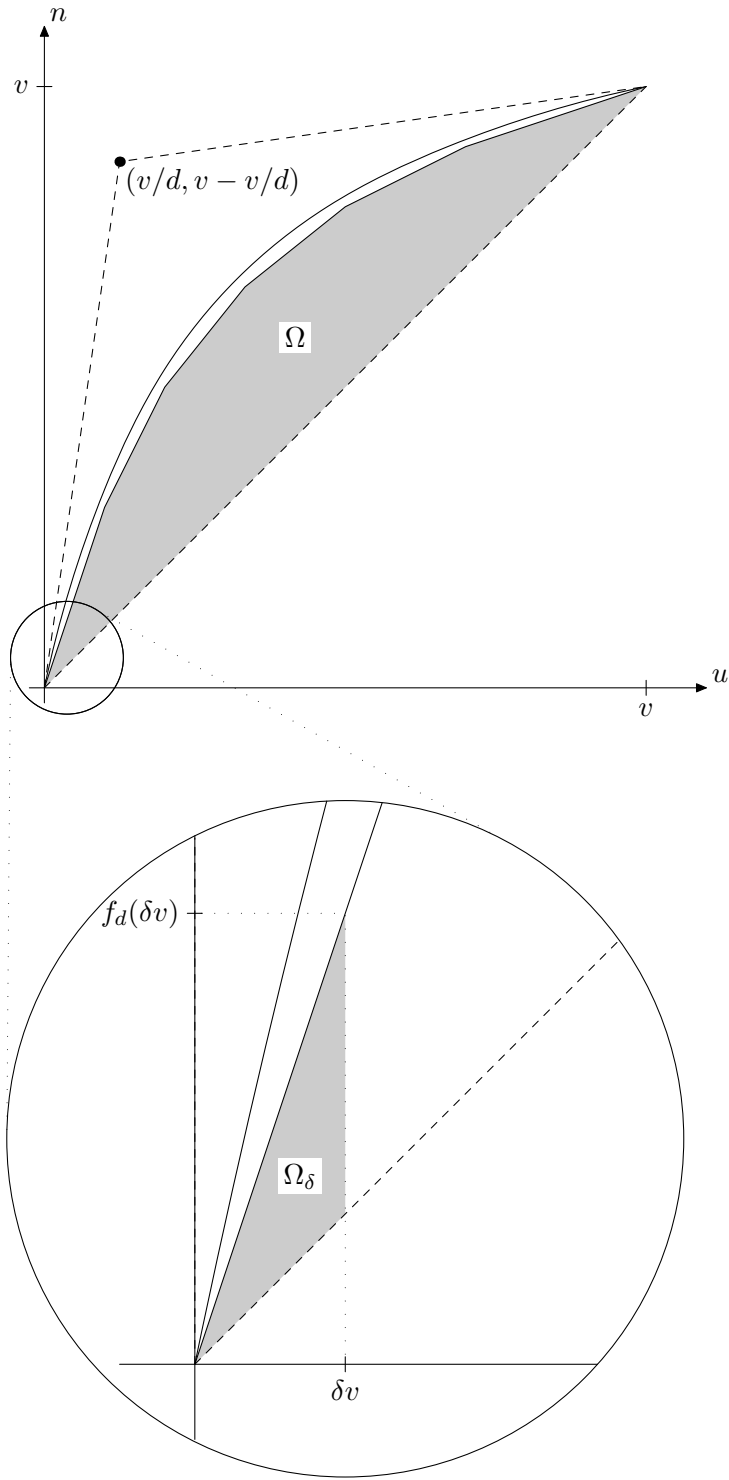}}
  \caption{The level curve $Q(u,n) = 1$ and the function $f_d$
    plotted for $d = 8$ with solid lines. The lower picture shows half
    of the set~$\Omega_{\delta}$.}
  \label{fig:rubbet}
\end{figure}

As can be seen from Figure~\ref{fig:rubbet}, $2|\Omega_{\delta}| \leq
2 \delta v f(\delta v)$ and $2|\Omega'| \leq v^2$. The proof now
proceeds by setting $\delta = 10^{-5}$ and then proving that
\begin{gather}
  2 \delta v f(\delta v) \max_{(u,n)\in\Omega_{\delta}} P(u,n) < 1/2
  \mbox{,}\\
  v^2 \max_{(u,n)\in\Omega'} P(u,n) < 1/2
  \mbox{.}
\end{gather}
This is enough to complete the proof, since in that case the
probability that a randomly selected graph has the desired properties
is non-zero and hence there exists at least one graph with the desired
properties.

\subsection{Analysis close to extreme points}
\label{sec:extreme}

Since both $P$ and our functions $f_d$ for $d \in \{5,6,7,8\}$ are
symmetric with respect to reflection around the line $u + n = v$, it
is enough to consider pairs $(u,n) \in \Omega_{\delta}$ such that $u +
n \leq v$. For fixed $u \leq \delta v$,
\begin{displaymath}
  \frac{P(u,n+1)}{P(u,n)}
  =
  \frac{(v-n)}{(n+1)}\prod_{i=1}^{d}\frac{dn+i}{dn-du+i}
  >
  1
  \mbox{,}
\end{displaymath}
therefore $P(u,n)$ is increasing in~$n$. Hence, it suffices to bound
$P(u,f_d(u))$. The following lemma establishes a slightly more general
result.
\begin{lemma}
  For every integer~$u$ such that $1 \leq u \leq 10^{-5}v$ and every
  $(k,d) \in \{(2,5),(3,6),(3,7),(3,8)\}$, $2k \delta^2 v^2 P(u,ku) <
  1/2$ where $P$~is defined by~(\ref{eq:P}).
\end{lemma}
\begin{proof}
  The proof is by induction on~$u$. The base case is clear since
  \begin{displaymath}
    2k \delta^2 v^2 P(1,k) = 
    2k(\delta v)^2 \binom{v}{1} \binom{v}{k} \frac{(kd)!(dv-d)!}{(dv)!(kd-d)!}
    < 2 \delta^2 \frac{k^{d+1}}{k!}\cdot\frac{v^{3+k}}{(v-1)^d}
    \,\mbox{.}
  \end{displaymath}
  For $d > 3+k$ and $v$~large enough, the latter expression above is
  strictly less than $1/2$. For the cases when $d = 3+k$, the latter
  expression is less than
  \begin{displaymath}
    10^{-10} \cdot 3^7 \cdot \frac{v^d}{(v-1)^d}
    <
    1/2
    \mbox{.}
  \end{displaymath}
  For the inductive step, we show that $P(u,ku) / P(u+1,ku+k) > 1$.
  Since
  \begin{displaymath}
    \frac{P(u,ku)}{P(u+1,ku+k)}
    =
    \frac{\binom{v}{u}\binom{v}{ku}\binom{kdu}{du}\big/\binom{dv}{du}}
    {\binom{v}{u+1}\binom{v}{ku+k}\binom{kdu+kd}{du+d}\big/\binom{dv}{du+d}}
  \end{displaymath}
  we need estimates of the following form:
  \begin{gather*}
    \binom{v}{u} \bigg/ \binom{v}{u+1} = \frac{u+1}{v-u} > \frac{u+1}{v}
    \,\mbox{,}\\
    \binom{v}{ku} \bigg/ \binom{v}{ku+k}
    = \frac{\prod_{i=1}^k (ku+i)}{\prod_{i=0}^{k-1}(v-ku-i)}
    > \frac{(ku+1)^k}{(v-ku)^k}
    > \frac{(u+1)^k}{v^k}
    \,\mbox{,}\\
    \binom{dv}{du+d} \bigg/ \binom{dv}{du}
    > \frac{(dv-du-d)^d}{(du+d)^d}
    > \frac{(v-2u)^d}{(u+1)^d}
    \,\mbox{,}\\
    \binom{kdu}{du} \bigg/ \binom{kdu+kd}{du+d}
    >
    \frac{(du)^{d}((k-1)du)^{(k-1)d}}{(kdu+kd)^{kd}}
    =
    \frac{(k-1)^{(k-1)d}u^{kd}}{k^{kd}(u+1)^{kd}}
    \,\mbox{.}
  \end{gather*}
  Put together, the above bounds imply that
  \begin{displaymath}
    \frac{P(u,ku)}{P(u+1,ku+k)}
    >
    \frac{(v-2u)^d}{v^{k+1}(u+1)^{d-k-1}}
    \cdot
    \frac{(k-1)^{(k-1)d}}{k^{kd}}
    \cdot
    \frac{u^{kd}}{(u+1)^{kd}}
  \end{displaymath}
  For $v > \delta^2$ and $u$ such that $1 \leq u \leq \delta v$, this
  is greater than
  \begin{multline*}
    \frac{(1-2\delta)^d}{(\delta+\delta^2)^{d-k-1}}
    \cdot
    \frac{(k-1)^{(k-1)d}}{k^{kd}}
    \cdot
    \frac{\delta^{kd}}{(\delta+\delta^2)^{kd}}
    \\
    >
    \delta^{-(d-k-1)}\bigl(1-(kd+3d)\delta\bigr)
    \frac{(k-1)^{(k-1)d}}{k^{kd}}
    \,\mbox{.}
  \end{multline*}
  For $(k,d) = (2,5)$ the above ratio is at least $10^{10} 2^{-10}
  (1-25\delta) > 1$.  For $k = 3$ the ratio is at least $10^{5d-20}
  2^{2d} 3^{-3d} (1-6d\delta)$, which is strictly greater than one
  for the $d$~considered.
\end{proof}

\subsection{The interior region}
\label{sec:interior}

To bound $v^2 P(u,n)$ in~$\Omega'$, write $u = \alpha v$ and $n =
\beta v$ and apply Stirling's formula
\begin{displaymath}
  \binom{v}{\alpha v} =
  \bigl(\alpha^{-\alpha} (1-\alpha)^{-(1-\alpha)}\bigr)^v
  \poly(v)
\end{displaymath}
to $P(u,n)$ defined in (\ref{eq:P}):
\begin{displaymath}
  v^2 P(\alpha v, \beta v) =
  \biggl(
  \frac{(1-\alpha)^{(d-1)(1-\alpha)} \beta^{(d-1)\beta}}
  {\alpha^{\alpha} (1-\beta)^{(1-\beta)}
    (\beta-\alpha)^{d(\beta-\alpha)}}
  \biggr)^v \poly(v)
  \mbox{.}
\end{displaymath}
Note that the above expression is valid also for $u = n$, i.e., also
for $\alpha = \beta$, if we use the convention that $0^0 = 1$. By the
symmetry of the function~$P$ and the results obtained in
\S~\ref{sec:extreme} it is enough to consider pairs $(\alpha, \beta)$
in the set
\begin{equation}
  \label{eq:A}
  A = \{(\alpha,\beta) \suchthat
  \beta-\alpha \geq 0 \land \beta v \leq f_d(\alpha v)
  \land \alpha \geq 10^{-5} \land \beta \leq 1-10^{-5}\}
  \mbox{.}
\end{equation}
Hence, it is sufficient to prove that there exists a universal
constant~$c < 1$, strictly bounded away from~$1$, such that
\begin{equation}
  \label{eq:Q}
  Q(\alpha, \beta) =
  \frac{(1-\alpha)^{(d-1)(1-\alpha)} \beta^{(d-1)\beta}}
       {\alpha^{\alpha} (1-\beta)^{(1-\beta)}
        (\beta-\alpha)^{d(\beta-\alpha)}}
  <
  c
  \quad\mbox{for $(\alpha,\beta) \in A$.}
\end{equation}
To this end, we first show that it is enough to consider the
boundaries of~$A$, i.e., the points where $\beta v = f_d(\alpha v)$,
and then analyze the function~$Q$ on those line segments.

\subsubsection{It is enough to consider the boundary}

To achieve the first goal, we prove that $\ln Q$ is convex along lines
of the form $\beta - \alpha = y$ for non-negative~$y$. This amounts to
substituting, e.g., $\alpha = (1+x-y)/2$ and $\beta = (1+x+y)/2$ in
the expression for~$Q$ and then consider the resulting expression as a
function of~$x$ for arbitrary fixed~$y$.
\begin{lemma}
  \label{lem:convex}
  For every fixed $y \in [0, 1-2/d]$, the function
  \begin{displaymath}
    q(x; y) \from x \mapsto \ln Q((1+x-y)/2, (1+x+y)/2)    
  \end{displaymath}
  is convex in the interval
  \begin{displaymath}
    |x| \leq \frac{d(1-y)-2}{d-2}
    \,\mbox{.}
  \end{displaymath}
\end{lemma}
\begin{proof}
  Straightforward substitution shows that $q(x; y) = -dy(\ln2 + \ln y)
  - (d-2)\ln2 + \frac{1}{2}g(x; y)$ where
  \begin{multline*}
    g(x; y) = (d-1)(1+x+y)\ln(1+x+y) + (d-1)(1-x+y)\ln(1-x+y)\\
    {} - (1+x-y)\ln(1+x-y) - (1-x-y)\ln(1-x-y)
    \mbox{.}
  \end{multline*}
  Hence the derivative of~$q$ with respect to~$x$ is
  \begin{displaymath}
    g'(x; y) = (d-1)\ln\frac{1+x+y}{1-x+y} - \ln\frac{1+x-y}{1-x-y}
  \end{displaymath}
  and the second derivative is, consequently,
  \begin{displaymath}
    g''(x; y)
    =
    \frac{(d-1)(2+2y)}{(1+x+y)(1-x+y)} - \frac{2-2y}{(1-x-y)(1+x-y)}
    \,\mbox{.}
  \end{displaymath}
  We now rewrite the second derivative as
  \begin{displaymath}
    g''(x; y)
    =
    2\frac{(d-2)(1-x^2-y^2) - dy(1+x^2-y^2)}{(1+x+y)(1-x+y)(1-x-y)(1+x-y)}
  \end{displaymath}
  and obtain that $g''(x; y)$ is non-negative when
  \begin{displaymath}
    (d-2)(1-x^2-y^2) \geq dy(1+x^2-y^2)
    \mbox{,}
  \end{displaymath}
  or, equivalently,
  \begin{displaymath}
    x^2(d(1+y)-2) \leq (d(1-y)-2)(1-y^2)
    \mbox{.}
  \end{displaymath}
  It is now straightforward to see by substitution that the above
  inequality is satisfied as soon as $|x| \leq (d(1-y)-2)/(d-2)$ and
  $y \in [0, 1 - 2/d]$. Finally, we remark that the same result is
  valid also for the case when $y = 0$: $q(x; 0)$ is convex for $|x|
  \leq 1$.
\end{proof}
To summarize, the function $q$ considered in the above lemma is convex
along lines parallel to the $x$-axis inside a triangle with corners
$(x,y) \in \{(1,0), (-1,0), (0,1-2/d)\}$. Translated to
$(\alpha,\beta)$-coordinates, the lemma therefore implies that it is
enough to bound the function $Q(\alpha,\beta)$ on the boundaries
of~$A$ as soon as $\Omega'$ is contained inside the triangle with
corners $(\alpha,\beta) \in \{(1,1),(0,0),(1/d,1-1/d)\}$. For our
case, it can be seen without much ado that $\Omega'$ is indeed
contained in this triangle: the functions $f_d$ are piecewise linear
and the slopes of the line segments are all strictly less than the
slope of the line from $(0,0)$ to $(1/d,1-1/d)$.

\subsubsection{Bounding the function on the boundary}

Since the function $Q(\alpha,\beta)$ is symmetric with respect to
reflection around the line $\alpha + \beta = 1$, it is sufficient to
prove that $Q(\alpha,\beta)$ is strictly less than one on the first
three ``legs'' of~$f_5$, the first three ``legs'' of~$f_6$, the first
four ``legs'' of~$f_7$, and the first three ``legs'' of~$f_8$. In
principle, this can be done by substituting $\beta = f_d(\alpha)$ and
then analyzing the resulting function using calculus. Since this is
extraordinarily tedious---both for the author and for the reader---we
instead choose a computer-assisted method of proof using interval
arithmetic~\cite{Moo66}. Specifically, we present a computer program
that verifies the following claims:
\begin{claim}
  \label{claim:5}
  For $d = 5$ and $\alpha \in [10^{-5},0.5]$, $Q(\alpha, \frac{1}{v}
  f_5(\alpha v)) \leq 0.9999$ where $Q$ is defined by~(\ref{eq:Q}) and
  $f_5$ by~(\ref{eq:f5}).
\end{claim}
\begin{claim}
  \label{claim:6}
  For $d = 6$ and $\alpha \in [10^{-5},0.5]$, $Q(\alpha, \frac{1}{v}
  f_6(\alpha v)) \leq 0.9999$ where $Q$ is defined by~(\ref{eq:Q}) and
  $f_6$ by~(\ref{eq:f6}).
\end{claim}
\begin{claim}
  \label{claim:7}
  For $d = 7$ and $\alpha \in [10^{-5},0.39]$, $Q(\alpha, \frac{1}{v}
  f_7(\alpha v)) \leq 0.9999$ where $Q$ is defined by~(\ref{eq:Q}) and
  $f_7$ by~(\ref{eq:f7}).
\end{claim}
\begin{claim}
  \label{claim:8}
  For $d = 8$ and $\alpha \in [10^{-5},0.34]$, $Q(\alpha, \frac{1}{v}
  f_8(\alpha v)) \leq 0.9999$ where $Q$ is defined by~(\ref{eq:Q}) and
  $f_8$ by~(\ref{eq:f8}).
\end{claim}
There is built-in support for interval arithmetic in some compilers
for some programming languages. The author has constructed a C++
program that uses the interval arithmetic routines built into Sun One
Studio~7~\cite{SunIntervalC++}. Complete source code of the program is
given in the appendix. Most of it is self-explanatory, the exception
being, maybe, the recursive function that verifies that some function
$f(x,y)$ is strictly less than some given upper bound~$b$ along the
line $y=y(x)$:
\begin{Verbatim}[fontsize=\footnotesize,xleftmargin=\leftmargini]
bool check(const AffineFunction& y, const di& i) {
    if(clt(f(i,y(i)), b)) {
        log << i << ": " << sup(f(i,y(i))) << std::endl;
        return true;
    }
    else {
        return check(y, interval_hull(di(inf(i)), di(mid(i)))) &&
               check(y, interval_hull(di(mid(i)), di(sup(i))));
    }
}
\end{Verbatim}
This function computes an interval~$I$ containing $\{f(x,y(x))
\suchthat x \in i\}$. If that interval is certainly less than an
interval~$b$ that contains the desired upper bound, i.e., if
$\forall z \in I \forall w \in b \,[z < w]$,
the claim has been verified on the interval~$i$. Otherwise, the
verification proceeds recursively: Two intervals that together
cover~$i$ are constructed, and the verification continues
on those intervals. A transcript of the verified subintervals
is written to a log file.

When compiled with Sun One Studio~7 and executed on a Sun workstation,
the program verified Claims~\ref{claim:5}--\ref{claim:8} in a couple
of seconds. The number of subintervals used were 332, 391, 857, and
261, respectively. The only option given to the compiler was
``\mbox{-xia}'', which enables support for interval arithmetic.

\section{Conclusions}

The methods described in this paper are fairly general and can be
applied to larger values of~$d$ without any complications. For
smaller~$d$, the analysis of the function~$P$ close to the extreme
points needs to be adapted.

Our functions~$f_d$ are all symmetric along the line $u + n = v$. This
property follows from our use of the probabilistic method; the
involved probabilities have the same symmetry. One possible direction
for future work could be to improve the behaviour of the
functions~$f_d$ in the region where $u$ is close to~$v$. Indeed, for
large~$d$, one would expect that the neighbour set of any $0.99v$
vertices is the entire other side of the bipartite graph.

\section{Acknowledgments}

Per-Olof Persson suggested, with great insight, to the author that he
should write the program verifying Claims~\ref{claim:5}--\ref{claim:8}
in C++ instead of Fortran. In addition, Staffan Gustafsson gave many
helpful comments that made the C++ program more readable.

\clearpage

\appendix

\section{Source code}

\begin{Verbatim}[fontsize=\footnotesize,numbers=left]
#include <fstream>
#include <suninterval.h>
typedef SUNW_interval::interval<double> di;

class Q
{
    const di d;
    const di one;
public:
    Q(const int degree) : d(degree), one("[1]") {}
    di operator()(const di& alpha, const di& beta) const {
        return (pow(one-alpha, (d-one)*(one-alpha)) * pow(beta,(d-one)*beta))
               / (pow(alpha,alpha) * pow(one-beta,one-beta)
                  * pow(beta-alpha, d*(beta-alpha)));
    }
};

class AffineFunction
{
    const di k;
    const di m;
public:
    AffineFunction(const di& slope, const di& offset)
        : k(slope), m(offset) {}
    AffineFunction(const di& x0, const di& y0, const di& x1, const di& y1)
        : k((y1-y0)/(x1-x0)), m(y0-k*x0) {}
    di operator()(const di& x) const { return k*x+m; }
};

class Segment
{
    const di i;
    const AffineFunction f;
public:
    Segment(const di& preimage, const di& slope, const di& offset)
        : i(preimage), f(slope, offset) {}
    Segment(const di& x0, const di& y0, const di& x1, const di& y1)
        : i(interval_hull(x0,x1)), f(x0, y0, x1, y1) {}
    const di& preimage() const { return i; }
    const AffineFunction& function() const { return f; }
};

class Checker
{
    const di b;
    const Q& f;
    std::ofstream log;
    bool check(const AffineFunction& y, const di& i) {
        if(clt(f(i,y(i)), b)) {
            log << i << ": " << sup(f(i,y(i))) << std::endl;
            return true;
        }
        else {
            return check(y, interval_hull(di(inf(i)), di(mid(i)))) &&
                   check(y, interval_hull(di(mid(i)), di(sup(i))));
        }
    }
public:
    Checker(const di& bound, const Q& fun, const char *file)
        : b(bound), f(fun), log(file) {}
    bool operator()(const Segment& s) {
        return check(s.function(), s.preimage());
    }
};

int main()
{
    const di bound("[0.9999]");
    // Verify claim for d==5
    Checker C5(bound, Q(5), "Q5.txt");
    Segment L51(di("[1e-5,0.15]"),  di("[2]"), di("[0]"));
    Segment L52(di("[0.15]"), di("[0.30]"),  di("[0.30]"), di("[0.50]"));
    Segment L53(di("[0.30]"), di("[0.50]"),  di("[0.50]"), di("[0.70]"));
    if(C5(L51) && C5(L52) && C5(L53)) {
        std::cout << "Claim is true for d==5." << std::endl;
    }
    // Verify claim for d==6
    Checker C6(bound, Q(6), "Q6.txt");
    Segment L61(di("[1e-5,0.10]"),  di("[2.5]"), di("[0]"));
    Segment L62(di("[0.10]"), di("[0.25]"),  di("[0.25]"), di("[0.50]"));
    Segment L63(di("[0.25]"), di("[0.50]"),  di("[0.50]"), di("[0.75]"));
    if(C6(L61) && C6(L62) && C6(L63)) {
        std::cout << "Claim is true for d==6." << std::endl;
    }
    // Verify claim for d==7
    Checker C7(bound, Q(7), "Q7.txt");
    Segment L71(di("[1e-5,0.10]"),  di("[3]"), di("[0]"));
    Segment L72(di("[0.10]"), di("[0.30]"),  di("[0.15]"), di("[0.40]"));
    Segment L73(di("[0.15]"), di("[0.40]"),  di("[0.30]"), di("[0.61]"));
    Segment L74(di("[0.30]"), di("[0.61]"),  di("[0.39]"), di("[0.70]"));
    if(C7(L71) && C7(L72) && C7(L73) && C7(L74)) {
        std::cout << "Claim is true for d==7." << std::endl;
    }
    // Verify claim for d==8
    Checker C8(bound, Q(8), "Q8.txt");
    Segment L81(di("[1e-5,0.10]"),  di("[3]"), di("[0]"));
    Segment L82(di("[0.10]"), di("[0.30]"),  di("[0.20]"), di("[0.50]"));
    Segment L83(di("[0.20,0.34]"),  di("[1.25]"), di("[0.25]"));
    if(C8(L81) && C8(L82) && C8(L83)) {
        std::cout << "Claim is true for d==8." << std::endl;
    }
}
\end{Verbatim}

\end{document}